\documentclass[11pt,leqno]{article}

\usepackage[french]{babel}
\usepackage[applemac]{inputenc}
\usepackage{amssymb,amsmath,amsthm,latexsym}

\numberwithin{equation}{section}
 
\let \a = \alpha
\def \d {\displaystyle}
\let \m = \medbreak
\let \n = \noindent
\def \Ind {\mathop {\Ind}\limits}

\def\build#1_#2^#3{\mathrel{\mathop{\kern 0pt#1}\limits_{#2}^{#3}}}

\def\Ind{\mathop{\hbox{Ind}}}

\def\smallcomp{\hbox{\fiverm C}\kern-.35em{\hbox{\fiverm I}}}
\def\smallreel{\hbox{\fiverm R}\kern-.6em{\hbox{\fiverm I}}}
\def\tore{\hbox{\rm T}\kern-.65em{\hbox{\rm I} } }
\def\fleche{\hbox{\sl v}\kern-.2em{\rightarrow}}
\def\equivalent{\Leftarrow\kern-.2em\Rightarrow}
\def\dess #1 by #2 (#3){
  \vbox to #2{
    \hrule width #1 height 0pt depth 0pt
    \vfill
    \special{picture #3} 
    }
  }

\def\dessin #1 by #2 (#3 scaled #4){{
  \dimen0=#1 \dimen1=#2
  \divide\dimen0 by 1000 \multiply\dimen0 by #4
  \divide\dimen1 by 1000 \multiply\dimen1 by #4
  \dess \dimen0 by \dimen1 (#3 scaled #4)}
  }

\def\hfl#1#2{\smash{\mathop{\hbox to
12mm{\rightarrowfill}}\limits^{\scriptstyle#1}_{\scriptstyle#2}}}

\def\nbZ{{\mathchoice {\hbox{$\sf\textstyle Z\kern-0.4em Z$}}
{\hbox{$\sf\textstyle Z\kern-0.4em Z$}}
{\hbox{$\sf\scriptstyle Z\kern-0.3em Z$}}
{\hbox{$\sf\scriptscriptstyle Z\kern-0.2em Z$}}}}

\def\nbQ{{\mathchoice {\setbox0=\hbox{$\displaystyle\rm
Q$}\hbox{\raise
0.15\ht0\hbox to0pt{\kern0.4\wd0\vrule height0.8\ht0\hss}\box0}}
{\setbox0=\hbox{$\textstyle\rm Q$}\hbox{\raise
0.15\ht0\hbox to0pt{\kern0.4\wd0\vrule height0.8\ht0\hss}\box0}}
{\setbox0=\hbox{$\scriptstyle\rm Q$}\hbox{\raise
0.15\ht0\hbox to0pt{\kern0.4\wd0\vrule height0.7\ht0\hss}\box0}}
{\setbox0=\hbox{$\scriptscriptstyle\rm Q$}\hbox{\raise
0.15\ht0\hbox to0pt{\kern0.4\wd0\vrule height0.7\ht0\hss}\box0}}}}

\def\nbC{{\mathchoice {\setbox0=\hbox{$\displaystyle\rm C$}%
\hbox{\hbox to0pt{\kern0.4\wd0\vrule height0.9\ht0\hss}\box0}}
{\setbox0=\hbox{$\textstyle\rm C$}\hbox{\hbox
to0pt{\kern0.4\wd0\vrule height0.9\ht0\hss}\box0}}
{\setbox0=\hbox{$\scriptstyle\rm C$}\hbox{\hbox
to0pt{\kern0.4\wd0\vrule height0.9\ht0\hss}\box0}}
{\setbox0=\hbox{$\scriptscriptstyle\rm C$}\hbox{\hbox
to0pt{\kern0.4\wd0\vrule height0.9\ht0\hss}\box0}}}}

\begin{document}

\centerline{\bf UNE NOUVELLE MINORATION DE LA TRACE DES ENTIERS ALGEBRIQUES}
\centerline{\bf TOTALEMENT POSITIFS}
 
\vskip1cm

\centerline { V. FLAMMANG }

\section{Introduction}
 \n
Soit $\alpha$   un entier algébrique de degré $d \geq 2$, totalement positif, i.e.,  dont les  conjugués  $\alpha=\alpha_1$, \ldots, $\alpha_d$ sont tous des réels positifs et de polynôme minimal $P$. On pose
$${\rm{S}}_k = \d \sum_{i=1}^{d} {\a}_i ^k.$$
\m
\n
La {\it trace absolue de $\alpha$} se définit par
$$\d {\cal T}race (\a )=\frac{1}{d}\text{trace}(\a)= \frac{1}{d} {\rm{S}}_1$$
\m
\n
et on désigne par $\cal T$ l'ensemble de tels $\d {\cal T}race (\a )$.\\
\m
\n
Le  {\it problème de Schur-Siegel-Smyth pour la trace} (appelé ainsi par P. Borwein dans son livre \cite{B}) est le suivant:\\

On se fixe $\rho < 2$. Il faut alors montrer que tous les entiers algébriques totalement\\
\indent positifs, sauf un nombre fini,  vérifient $\d {\cal T}race (\a )> \rho$.\\
\m
\n
Le problème a été résolu en 1918 par I. Schur pour $\rho < \sqrt e $ \cite{Sc} puis en 1945 par C. L. Siegel pour $\rho < 1.7337$ \cite{Si}. Les résultats de Schur et de Siegel font intervenir des inégalités sur le discriminant d'un entier algébrique qui est la quantité
$${\rm Disc}(\alpha)= \d \prod_{1 \leq i < j \leq d} (\alpha_i - \alpha_j)^2.$$
\m
\n
Tous les nombreux résultats ultérieurs se basent sur le principe des {\it fonctions auxiliaires} qui repose sur le fait que le résultant de deux polynômes à coefficients entiers et sans facteurs communs est un entier non nul. Ainsi en 1984, C.J. Smyth résout le problème pour $\rho < 1.7719$ \cite{Sm3}, en 1997, V. Flammang, M. Grandcolas et G. Rhin pour $\rho < 1.7735$ \cite{F2}, en 2004, J. McKee et C.J. Smyth pour $\rho < 1.7783786$ \cite{McS1}, en 2006, J. Aguirre, M. Bilbao et J. C. Peral pour $\rho < 1.7800$ \cite{ABP}, en 2006 par V. Flammang pour $\rho < 1.7822$ (communication privée à C. J. Smyth). En 2007, J. Aguirre et J. C. Peral ont résolu le problème pour $\rho< 1.7836$ \cite{AP1} puis en 2008 pour $\rho < 1.784109$ \cite{AP2}. En 2009, nous l'avons résolu pour $\rho < 1.78702$ \cite{F3}. J. McKee \cite{Mc} a résolu le problème en 2011 pour $\rho < 1.78839$ en utilisant une partie de nos polynômes ayant des racines complexes provenant de la fonction auxiliaire donnée dans \cite{F3}.  En 2011 également, Y. Liang et Q. Wu \cite{LW} ont résolu le problème pour $\rho < 1.79193$.\\
 Dans cette note, nous résolvons le problème pour $\rho < 1,792812$. Plus précisément, nous montrons ici le théorème suivant:

\vfill \eject

\newtheorem{guess}{Théorème}

\begin{guess}
Si $\a$ est un entier algébrique totalement positif de degré $d$ et polynôme minimal différent de $x$, $x-1$, $x^2 - 3x + 1$,  $x^3 - 5x^2 + 6x - 1$, $x^4 - 7x^3 + 13x^2 - 7x + 1$ et $x^4 - 7x^3 + 14x^2 - 8x + 1$ alors on a:
\end{guess}
$$\frac{1}{d}\text{trace}(\a) \geq 1.792812.$$

\vskip1cm

\n
Rappelons qu'un entier algébrique $\alpha$ de degré $d$ est {\it réciproque} si son polynôme minimal $P$ vérifie $\d P(x) =x^d. P(1/x)$.\\

\newtheorem{corollaire}{Corollaire}
\begin{corollaire}
Si $\a$ est un entier algébrique totalement positif et réciproque de degré $d$ et de polynôme minimal différent de $x^2 - 3x + 1$ et $x^4 - 7x^3 + 13x^2 - 7x + 1$, alors on a:
\end{corollaire}
$$\frac{1}{d}\text{trace}(\a) \geq 1.896406.$$

\vskip1cm
\n
{\bf Remarques}\\
\n
La minorations  da Théorème 1 améliore la plus récente connue à savoir celle obtenue par Y. Liang et Q. Wu  en 2011 qui est 1.79193.\\
\m
\n
La minoration du corollaire améliore celle de X. Dong et Q. Wu \cite{DW} qui est de 1.8945909.
\m
\n
Les polynômes qui interviennent dans ce théorème ainsi que leurs exposants se trouvent en fin de note.

 \section{Le principe des fonctions auxiliaires}
 
\m
\n
La  fonction auxiliaire qui intervient dans l'étude de la trace est du type:
 
$$ \mbox{pour}{\ } x > 0, {\ }f(x)\\= x - \sum_{ 1\leq j \leq J} c_j \log | Q_j(x)\\ |    {\ \ \ }(1)\\$$
 \n
où les $c_j$ sont des nombres réels positifs et les polynômes $Q_j$ sont des polynômes non nuls de $\nbZ[x]$.
\m
\n
Cette fonction auxiliaire a été introduite par C. J. Smyth dans \cite{Sm2}.\\
\m
\n Soit  $m$  le minimum de la fonction $f$. Si $P$ ne divise aucun des polynômes $Q_j$, on a alors
$$ \d \sum_{i=1}^{d} f( {\a}_i)\\ \geq md$$
\n
i.e.,
$$\text{trace}(\a)\\ \geq md + \sum_{ 1\leq j \leq J} c_j \log | \prod_{i=1}^{d} Q_j( {\a}_i)\\ | .$$

\vfill \eject
\n
Puisque $P$ ne divise aucun des  $Q_j$, alors $\d \prod_{i=1}^{d} Q_j( {\a}_i)$ est un entier non nul  car il est le résultant de $P$ et de $Q_j$.\\

\n
Par conséquent,  si $\a$ n'est pas une racine de $Q_j$, on a
$${\cal T}race( \a )\\ \geq m.  $$
\m
\n
Par ailleurs, J. P. Serre ( voir Appendix B dans \cite{AP2}) a montré que cette méthode ne peut pas donner une telle inégalité pour $\rho$ plus grand que 1.8983021\ldots Par conséquent, cette méthode ne peut pas être utilisée pour montrer que 2 est le plus petit point d'accumulation de $\cal T$.\\
\n
 Néanmoins, il est intéressant d'essayer d'obtenir des minorations pour ${\cal T}race (\a )$. Par exemple, cela a été utilisé dans la recherche de nombres de Salem de petit degré et de trace -2 par J. McKee et C. J. Smyth \cite{McS1}. Y. Liang et Q. Wu ont obtenu, grâce à leur résultat, le degré minimum d'un nombre de Salem de trace -4 et -5 \cite{LW}.
 
\section{Relation entre fonction auxiliaire et diamètre transfini entier généralisé}
\subsection{Rappels sur le diamètre transfini entier généralisé}
\n
Soit $K$ un compact de $\mathbb{C}$. Le {\it diamètre transfini de K} se définit par
$$\begin{array}{ccccc}
t(K)& = & \liminf  & \inf & |P |_{ \infty, K} ^{\frac{1}{n}}\\
                                              &     &   n \geq 1&  P \in \mathbb{C}[X] \\
                                              &    &  n \rightarrow   \infty & P{\  } \textsf{unitaire} \\
                                              &   &                                         & \deg(P)=n  
 \end{array} $$
où $|P |_{ \infty, K} = \d \sup_{z \in K} |P(z) |$ pour $P \in \mathbb{C} [X]$.\\
\m
\n
Nous définissons le {\it diamètre transfini entier  de K}  par
$$\begin{array}{ccccc} 
 t_{\nbZ}(K)& = & \liminf  &\inf & | P |_{\infty, K}^{\frac{1}{n}}              {\ \ \ \ \ } \\
                                              &     &   n \geq 1&  P \in \nbZ[X] &\\
                                              &    &  n \rightarrow   \infty & \deg(P)=n &
 \end{array} $$

\m
\n
Enfin, si $\varphi$ est une fonction positive définie sur $K$, le {\it $\varphi$-diamètre transfini entier généralisé de $K$} se définit par
$$\begin{array}{cccccc}
 t_{\mathbb{Z}, \varphi}(K)& = & \liminf  & \inf & \sup&\left ( | P(z) |^{\frac{1}{n}} {\ } \varphi(z) \right ).\\
                                              &     &   n \geq 1&  P \in \mathbb{Z}[X] & z \in K &\\
                                              &    &  n \rightarrow   \infty & \deg(P)=n &
 \end{array} $$
\m
\n
Cette version   du diamètre transfini entier pondéré a été introduite par F. Amoroso \cite{A2} et est un outil important dans l'étude des approximations rationnelles de logarithmes de nombres rationnels.

\vfill \eject

\subsection{Lien avec les fonctions auxiliaires}

\n
Dans la fonction auxiliaire (1), on remplace les $c_j$ par des nombres rationnels $a_j/q$ où $q$ est un entier  $>0$ tel que $q.c_j$ soit un entier pour tout $1\leq j \leq J$. On peut alors écrire:
$$ \mbox{pour}{\ } x > 0, {\ }f(x)\\= x -  \frac{t}{r} \log | Q(x)\\ | \geq m   {\ \ \ }(2)\\$$
\n
où $\d Q=\prod_{j=1}^J Q_j^{a_j} \in \nbZ[X]$ est de degré $\d r = \sum_{j=1}^J a_j \deg Q_j$ et $\d t= \sum_{j=1}^J c_j \deg Q_j$ (cette formulation a été introduite par J. P. Serre). Par conséquent, on cherche un polynôme $Q \in \nbZ[X]$ tel que
$$\sup_{x>0}~\vert Q(x) \vert^{t/r} e^{-x} \leq e^{-m}.$$
 
\n
Si l'on suppose t fixé, cela revient à trouver une borne supérieure effective pour le diamètre transfini entier pondéré de l'intervalle $[0, \infty[$ avec le poids $\varphi(x) = e^{-x}$:
$$\begin{array}{cccccc}
 t_{\nbZ,\varphi}([0,\infty [)& = & \liminf  & \inf & \sup&\left ( | P(x) |^{\frac{t}{r}} {\ }\varphi(x) \right )\\
                                              &     &   r \geq 1&  P \in \nbZ[X] & x>0 &\\
                                              &    &  r \rightarrow   \infty & \deg(P)=r &
 \end{array} $$
 \m
 \n
 {\bf Remarque:} Même si l'on a remplacé le compact $K$ par l'intervalle infini $[0, \infty [$, le poids $\varphi$ assure que la quantité  $t_{\nbZ,\varphi}([0, \infty[ )$ est finie.

\section{Construction d'une fonction auxiliaire}
\n
Le point essentiel est   de trouver une liste de  "bons" polynômes $Q_j$, i.e., qui donnent la meilleure valeur possible de $m$.  Jusqu'en 2003, les polynômes étaient trouvés de façon heuristique. Par exemple, dans \cite{Sm3} et \cite{AP1}, les auteurs ont cherché une collection de polynômes de petite trace absolue dont toutes les racines sont positives. 
\n
En 2003, Q. Wu \cite{Wu} a mis au point un algorithme qui permet une recherche systématique des "bons" polynômes. La méthode était la suivante. On considère une fonction auxiliaire comme celle définie en (1). On se fixe un ensemble $E_0$ de points de contrôle uniformément répartis sur un intervalle réel $I=[0; A]$ où $A$ est " suffisamment grand". Par LLL, on trouve un polynôme $Q$ petit sur $E_0$ au sens de la norme quadratique. On teste ce polynôme dans la fonction auxiliaire et on ne conserve que les facteurs de $Q$ qui ont un exposant non nul. La convergence de cette nouvelle fonction fournit des minima locaux que l'on ajoute à l'ensemble de points $E_0$ pour obtenir un nouvel ensemble de points de contrôle $E_1$. On relance LLL avec l'ensemble $E_1$ et on réitère le procédé.\\
\m
\n
En 2006, nous avons apporté deux améliorations à l'algorithme précédent dans l'utilisation de LLL. La première consiste, à chaque pas, à prendre en compte non seulement les nouveaux points de contrôle mais également les nouveaux polynômes de la fonction auxiliaire qui est la meilleure. La seconde est l'introduction d'un coefficient correcteur $t$.  L'idée est d'obtenir les bons polynômes $Q_j$ par récurrence.  Ainsi,  nous appelons cet algorithme {\it l'algorithme récursif}. Nous le détaillons, toujours pour la trace. Le premier pas consiste à optimiser la fonction  auxiliaire $f_1= x - t \log x$. On a alors $t=c_1$ où $c_1$  est la valeur qui donne la meilleure fonction $f_1$. On suppose qu'on a des polynômes $Q_1$, $Q_2$, \ldots, $Q_J$ et une fonction $f$ la meilleure possible pour cet ensemble de polynômes sous la forme (2). On cherche un polynôme $R \in \nbZ[x]$ de degré $k$ ($k=10$ par exemple) tel que
$$ \sup_{x \in I} | Q(x) R(x) | ^{\frac{t}{r+k}} e^{-x} \leq e^{-m},$$
\vfill \eject
\n
où $\d Q= \prod_{j=1}^J Q_j$.
\n
On veut donc que

$$  \sup_{x \in I} | Q(x) R(x) | \exp \left( \frac{-x(r+k)}{t} \right)$$
\n
soit aussi petit que possible. On applique LLL aux formes linéaires
$$ Q(x_i) R(x_i) \exp \left( \frac{-x_i(r+k)}{t} \right).$$
\n
Les $x_i$ sont des points de contrôle constitués de points uniformément répartis sur l'intervalle $I$ auxquels on a ajouté les points où $f$ a des minima locaux. On trouve donc un polynôme $R$ dont les facteurs irréductibles $R_j$ sont de bons candidats pour agrandir l'ensemble $\{Q_1, \ldots, Q_J\}$. On ne conserve que les facteurs $R_j$ qui ont un coefficient non nul dans la nouvelle fonction auxiliaire optimisée $f$. Après optimisation, certains  polynômes  précédents $Q_j$ peuvent avoir un exposant nul et sont alors rejetés.

\n
En 2009, pour obtenir la constante 1.78702, l'algorithme précédent avait été répété pour $k$ variant de 10 à 30. Nous nous étions arrêtés quand deux pas consécutifs ne fournissaient plus aucun nouveau polynôme. Ici, nous avons procédé différemment. Nous avons fait varier $k$ systématiquement de 3 à 82 et nous avons constaté que, même après plusieurs pas infructueux, des polynômes à exposant non nul apparaissent encore. C'est cette légère variante qui a permis l'amélioration des constantes.
\m
\n
{\bf Remarque}\\
\n
La constante $t$ obtenue pour la fonction auxiliaire du Théorème 1 vaut 2.6417021. Elle est proche de la valeur optimale de J. P. Serre qui est égale à 2.249214.\\

\section{Optimisation des $c_j$}
\n
On a à résoudre un problème du type suivant: trouver
$$ \max_{C} \min_{x \in X} f(x,C)$$
\n
où $f(x,C)$ est une fonction linéaire par rapport à $C=(c_0, c_1,\ldots, c_k)$ ($c_0$ est le coefficient de $x$ qui est égal à 1) et $X$ est un domaine compact de $\mathbb{C}$, le maximum étant pris pour $c_j \geq 0$ pour $j=0,\ldots, k$.
\n
Une solution classique consiste à prendre de très nombreux points de contrôle $(x_i)_{1 \leq i \leq N} $ et à résoudre le problème de programmation linéaire standard:
$$ \max_{C} \min_{1 \leq j \leq N} f(x_i,C).$$
 
 \n
 Mais alors le résultat obtenu dépend des points de contrôle choisis.\\
 \m
 \n
 L'idée de la programmation linéaire linéaire semi-infinie (introduite en Théorie des Nombres par C. J. Smyth \cite{Sm2}) consiste à répéter le processus précédent en ajoutant à chaque étape de nouveaux points de contrôle et à vérifier que ce procédé converge vers $m$, la valeur de la forme linéaire pour un choix de $C$ optimum. L'algorithme est le suivant:
 
(1) On choisit une valeur initiale de $C$ soit $C^0$ et on calcule
$$ m'_0= \min_{x \in X} f(x, C^0).$$

\vfill \eject
(2) On choisit un ensemble fini  $X_0$ de points de contrôle appartenant à X et l'on a\\
 $$m'_0 \leq m \leq m_0= \min_{x \in X_0}f(x, C^0).$$
 
 (3) On ajoute à $X_0$ les points où $f(x,C^0)$ admet des minima locaux pour obtenir un \\
 \hspace*{0.45in} nouvel ensemble $X_1$ de points de contrôle.
 
 (4) On résout le problème de programmation linéaire standard:
 $$ \max_{C} \min_{x \in X_1} f(x,C)$$
 
 \n 
On obtient une nouvelle valeur de $C$ notée $C^1$ et un résultat de la programmation linéaire égal à $m'_1=\d \min_{x \in X} f(x, C^1)$. On a alors
  $$m'_0 \leq m'_1 \leq m \leq m_1=\d \min_{x \in X_1}f(x, C^1) \leq m_0,$$
  
  (5) On répète les étapes (2) à (4) et on obtient donc deux suites $(m_i)$ et $(m'_i)$ qui \\
 \hspace*{0.45in}  vérifient
  $$ m'_0 \leq m'_1 \leq \ldots \leq m'_i \leq m \leq m_i \leq \ldots \leq m_1 \leq m_0,$$
  
  \n
 On s'arrête dès qu'il y a assez bonne convergence, quand par exemple, $m_i-m'_i \leq 10^{-6}$.\\
 \n
Supposons que $p$ itérations suffisent alors on prend $m=m'_p$.

\section{Preuve du Corollaire 1}
\m
\n
Pour obtenir leur résultat, X. Dong et Q. Wu \cite{DW} ont utilisé une fonction auxiliaire du type:
$$ \mbox{pour}{\ } x > 0, {\ }f(x)\\= \d x + \frac{1}{x} - \sum_{ 1\leq j \leq J} c_j \log | Q_j(x)| - \sum_{ 1\leq j \leq J} c_j \log | Q_j( \frac{1}{x})|. $$
\m
\n
Nous procédons différemment. Soit $\alpha$ un entier algébrique  totalement positif et réciproque de degré $d$ et de polynôme minimal $P$. Alors il existe un polynôme $Q$ totalement positif de degré $d/2$ vérifiant:\\
$$ P(X) = X^{d/2} Q( X + \frac{1}{X} -2).$$
\m
\n
Soient $\alpha_1$, $\ldots$, $\alpha_d$ et $\beta_1$, $\ldots$, $\beta_{d/2}$ les racines de $P$ et $Q$ respectivement. Alors on a:
$$  \mbox{pour}{\ }1 \leq i \leq d/2,{\ } \beta_i= \alpha_i + \frac{1}{\alpha_i} -2{\ \ \ }$$
\m
\n
Par conséquent, $\mbox{trace}(P)= \d \sum_{i=1}^d \alpha_i= \sum_{i=1}^{d/2} \left( \alpha_i + \frac{1}{\alpha_i} \right )= \sum_{i=1}^{d/2} (\beta_i + 2)$\\
\n
i.e.,  $\mbox{trace}(P)= \d \sum_{i=1}^{d/2} \beta_i + d = \mbox{trace}(Q) + d$. Or, le polynôme $Q$ vérifie les hypothèses du Théorème 1 donc\\
$$ \frac{\mbox{trace}(P)}{d} =  \frac{\mbox{trace}(Q)}{\frac{d}{2}.2} + 1 \geq \frac{1.792812}{2}+1= 1.896406.$$

\vfill \eject

\begin{center}
\bf{ Polynômes et exposants qui interviennent dans le Théorème 1}
\end{center}
\scriptsize
\n
 pol=$[x,\\
x - 1,\\
x - 2,\\
x^2 - 3 x + 1,\\
x^2 - 4 x + 1,\\
x^2 - 4 x + 2,\\
x^3 - 5 x^2 + 6 x - 1,\\
x^3 - 6 x^2 + 9 x - 3,\\
x^3 - 6 x^2 + 9 x - 1,\\
x^3 - 6 x^2 + 8 x - 1,\\
x^4 - 7 x^3 + 13 x^2 - 7 x + 1,\\
x^4 - 7 x^3 + 14 x^2 - 8 x + 1,\\
x^5 - 9 x^4 + 27 x^3 - 32 x^2 + 13 x - 1,\\
x^5 - 9 x^4 + 27 x^3 - 31 x^2 + 12 x - 1,\\
x^5 - 9 x^4 + 28 x^3 - 35 x^2 + 15 x - 1,\\
x^5 - 9 x^4 + 26 x^3 - 29 x^2 + 11 x - 1,\\
x^6 - 11 x^5 + 43 x^4 - 72 x^3 + 51 x^2 - 14 x + 1,\\
x^6 - 11 x^5 + 43 x^4 - 73 x^3 + 53 x^2 - 15 x + 1,\\
x^7 - 13 x^6 + 64 x^5 - 150 x^4 + 172 x^3 - 89 x^2 + 18 x - 1,\\
x^7 - 13 x^6 + 63 x^5 - 143 x^4 + 157 x^3 - 78 x^2 + 16 x - 1,\\
x^7 - 12 x^6 + 54 x^5 - 114 x^4 + 117 x^3 - 56 x^2 + 12 x - 1,\\
x^7 - 13 x^6 + 63 x^5 - 144 x^4 + 160 x^3 - 80 x^2 + 16 x - 1,\\
x^8 - 15 x^7 + 89 x^6 - 269 x^5 + 445 x^4 - 402 x^3 + 187 x^2 - 40 x + 3,\\
x^8 - 15 x^7 + 90 x^6 - 277 x^5 + 467 x^4 - 428 x^3 + 200 x^2 - 42 x + 3,\\
x^8 - 14 x^7 + 78 x^6 - 221 x^5 + 339 x^4 - 277 x^3 + 111 x^2 - 19 x + 1,\\
x^8 - 14 x^7 + 78 x^6 - 222 x^5 + 345 x^4 - 289 x^3 + 120 x^2 - 21 x + 1,\\
x^9 - 16 x^8 + 103 x^7 - 344 x^6 + 643 x^5 - 681 x^4 + 399 x^3 - 123 x^2 + 18 x - 1,\\
x ^{10} - 18 x^9 + 134 x^8 - 537 x^7 + 1265 x^6 - 1798 x^5 + 1526 x^4 - 743 x^3 + 194 x^2 - 24 x + 1,\\
x ^{10} - 18 x^9 + 134 x^8 - 538 x^7 + 1273 x^6 - 1822 x^5 + 1560 x^4 - 766 x^3 + 200 x^2 - 24 x + 1,\\
x ^{10} - 18 x^9 + 135 x^8 - 549 x^7 + 1320 x^6 - 1920 x^5 + 1662 x^4 - 813 x^3 + 206 x^2 - 24 x + 1,\\
x ^{10} - 19 x^9 + 150 x^8 - 643 x^7 + 1641 x^6 - 2573 x^5 + 2472 x^4 - 1412 x^3 + 451 x^2 - 71 x + 4,\\
x ^{10} - 18 x^9 + 135 x^8 - 549 x^7 + 1321 x^6 - 1929 x^5 + 1689 x^4 - 847 x^3 + 223 x^2 - 26 x + 1,\\
x ^{12} - 21 x ^{11} + 190 x ^{10} - 972 x^9 + 3103 x^8 - 6439 x^7 + 8780 x^6 - 7789 x^5 + 4372 x^4 - 1483 x^3 + 283 x^2 - 27 x + 1,\\
x ^{12} - 22 x ^{11} + 208 x ^{10} - 1108 x^9 + 3666 x^8 - 7840 x^7 + 10948 x^6 - 9877 x^5 + 5589 x^4 - 1885 x^3 + 349 x^2 - 31 x + 1,\\
x ^{12} - 21 x ^{11} + 190 x ^{10} - 972 x^9 + 3102 x^8 - 6430 x^7 + 8750 x^6 - 7742 x^5 + 4336 x^4 - 1470 x^3 + 281 x^2 - 27 x + 1,\\
x ^{12} - 22 x ^{11} + 207 x ^{10} - 1093 x^9 + 3574 x^8 - 7539 x^7 + 10373 x^6 - 9219 x^5 + 5145 x^4 - 1718 x^3 + 318 x^2 - 29 x + 1,\\
x ^{12} - 21 x ^{11} + 190 x ^{10} - 971 x^9 + 3090 x^8 - 6373 x^7 + 8613 x^6 - 7565 x^5 + 4216 x^4 - 1432 x^3 + 277 x^2 - 27 x + 1,\\
x ^{12} - 22 x ^{11} + 208 x ^{10} - 1108 x^9 + 3667 x^8 - 7851 x^7 + 10995 x^6 - 9977 x^5 + 5702 x^4 - 1952 x^3 + 368 x^2 - 33 x + 1,\\
x ^{12} - 22 x ^{11} + 209 x ^{10} - 1123 x^9 + 3757 x^8 - 8125 x^7 + 11435 x^6 - 10317 x^5 + 5775 x^4 - 1913 x^3 + 349 x^2 - 31 x + 1,\\
x ^{13} - 23 x ^{12} + 229 x ^{11} - 1298 x ^{10} + 4637 x^9 - 10930 x^8 + 17323 x^7 - 18505 x^6 + 13193 x^5 - 6143 x^4 + 1798 x^3 - 310 x^2 + 28 x - 1,\\
x ^{13} - 23 x ^{12} + 230 x ^{11} - 1313 x ^{10} + 4730 x^9 - 11240 x^8 + 17929 x^7 - 19217 x^6 + 13689 x^5 - 6338 x^4 + 1837 x^3 - 313 x^2 + 28 x - 1,\\
x ^{13} - 23 x ^{12} + 229 x ^{11} - 1298 x ^{10} + 4637 x^9 - 10930 x^8 + 17322 x^7 - 18499 x^6 + 13181 x^5 - 6134 x^4 + 1796 x^3 - 310 x^2 + 28 x - 1,\\
x ^{13} - 23 x ^{12} + 231 x ^{11} - 1333 x ^{10} + 4894 x^9 - 11963 x^8 + 19810 x^7 - 22195 x^6 + 16544 x^5 - 7942 x^4 + 2334 x^3 - 388 x^2 + 32 x - 1,\\
x ^{13} - 23 x ^{12} + 231 x ^{11} - 1333 x ^{10} + 4894 x^9 - 11962 x^8 + 19800 x^7 - 22157 x^6 + 16475 x^5 - 7880 x^4 + 2308 x^3 - 384 x^2 + 32 x - 1,\\
x ^{13} - 23 x ^{12} + 231 x ^{11} - 1333 x ^{10} + 4895 x^9 - 11975 x^8 + 19867 x^7 - 22332 x^6 + 16721 x^5 - 8062 x^4 + 2372 x^3 - 392 x^2 + 32 x - 1,\\
x ^{14} - 25 x ^{13} + 274 x ^{12} - 1736 x ^{11} + 7061 x ^{10} - 19368 x^9 + 36643 x^8 - 48110 x^7 + 43567 x^6 - 26760 x^5 + 10844 x^4 - 2778 x^3 + 420 x^2 - 33 x + 1,\\
x ^{14} - 25 x ^{13} + 274 x ^{12} - 1736 x ^{11} + 7060 x ^{10} - 19355 x^9 + 36573 x^8 - 47907 x^7 + 43221 x^6 - 26404 x^5 + 10625 x^4 - 2701 x^3 + 406 x^2 - 32 x + 1,\\
x ^{14} - 25 x ^{13} + 274 x ^{12} - 1736 x ^{11} + 7061 x ^{10} - 19369 x^9 + 36655 x^8 - 48167 x^7 + 43704 x^6 - 26937 x^5 + 10964 x^4 - 2816 x^3 + 424 x^2 - 33 x + 1,\\
x ^{14} - 25 x ^{13} + 275 x ^{12} - 1755 x ^{11} + 7216 x ^{10} - 20082 x^9 + 38697 x^8 - 51961 x^7 + 48328 x^6 - 30594 x^5 + 12779 x^4 - 3344 x^3 + 502 x^2 - 37 x + 1,\\
x ^{14} - 25 x ^{13} + 275 x ^{12} - 1754 x ^{11} + 7200 x ^{10} - 19973 x^9 + 38284 x^8 - 51007 x^7 + 46942 x^6 - 29332 x^5 + 12084 x^4 - 3131 x^3 + 472 x^2 - 36 x + 1,\\
x ^{14} - 25 x ^{13} + 274 x ^{12} - 1737 x ^{11} + 7077 x ^{10} - 19477 x^9 + 37056 x^8 - 49065 x^7 + 44959 x^6 - 28036 x^5 + 11556 x^4 - 3003 x^3 + 455 x^2 - 35 x + 1,\\
x ^{14} - 25 x ^{13} + 274 x ^{12} - 1735 x ^{11} + 7046 x ^{10} - 19275 x^9 + 36333 x^8 - 47505 x^7 + 42862 x^6 - 26281 x^5 + 10667 x^4 - 2747 x^3 + 418 x^2 - 33 x + 1,\\
x ^{14} - 25 x ^{13} + 274 x ^{12} - 1736 x ^{11} + 7061 x ^{10} - 19370 x^9 + 36666 x^8 - 48216 x^7 + 43819 x^6 - 27092 x^5 + 11086 x^4 - 2870 x^3 + 436 x^2 - 34 x + 1,\\
x ^{14} - 25 x ^{13} + 274 x ^{12} - 1735 x ^{11} + 7044 x ^{10} - 19248 x^9 + 36183 x^8 - 47058 x^7 + 42081 x^6 - 25460 x^5 + 10155 x^4 - 2568 x^3 + 387 x^2 - 31 x + 1,\\
x ^{14} - 25 x ^{13} + 274 x ^{12} - 1735 x ^{11} + 7044 x ^{10} - 19246 x^9 + 36160 x^8 - 46953 x^7 + 41834 x^6 - 25137 x^5 + 9920 x^4 - 2478 x^3 + 371 x^2 - 30 x + 1,\\
x ^{14} - 25 x ^{13} + 274 x ^{12} - 1735 x ^{11} + 7045 x ^{10} - 19262 x^9 + 36264 x^8 - 47311 x^7 + 42547 x^6 - 25976 x^5 + 10489 x^4 - 2685 x^3 + 406 x^2 - 32 x + 1,\\
x ^{15} - 27 x ^{14} + 322 x ^{13} - 2237 x ^{12} + 10058 x ^{11} - 30777 x ^{10} + 65693 x^9 - 98793 x^8 + 104703 x^7 - 77645 x^6 + 39727 x^5 - 13710 x^4 + 3079 x^3 - 424 x^2 + 32 x - 1,\\
x ^{15} - 27 x ^{14} + 323 x ^{13} - 2261 x ^{12} + 10304 x ^{11} - 32190 x ^{10} + 70713 x^9 - 110296 x^8 + 121940 x^7 - 94465 x^6 + 50220 x^5 - 17760 x^4 + 3993 x^3 - 532 x^2 + 37 x - 1,$
\vfill \eject
\n $
x ^{15} - 28 x ^{14} + 348 x ^{13} - 2533 x ^{12} + 11998 x ^{11} - 38873 x ^{10} + 88169 x^9 - 140945 x^8 + 157932 x^7 - 122084 x^6 + 63478 x^5 - 21465 x^4 + 4521 x^3 - 560 x^2 + 37 x - 1,\\
x ^{15} - 27 x ^{14} + 323 x ^{13} - 2261 x ^{12} + 10304 x ^{11} - 32189 x ^{10} + 70700 x^9 - 110227 x^8 + 121746 x^7 - 94152 x^6 + 49927 x^5 - 17607 x^4 + 3953 x^3 - 528 x^2 + 37 x - 1,\\
x ^{16} - 29 x ^{15} + 374 x ^{14} - 2836 x ^{13} + 14091 x ^{12} - 48408 x ^{11} + 118290 x ^{10} - 208431 x^9 + 265782 x^8 - 244226 x^7 + 159913 x^6 - 73242 x^5 + 22819 x^4 - 4635 x^3 + 573 x^2 - 38 x + 1,\\
x ^{16} - 29 x ^{15} + 374 x ^{14} - 2835 x ^{13} + 14069 x ^{12} - 48205 x ^{11} + 117251 x ^{10} - 205149 x^9 + 259051 x^8 - 235052 x^7 + 151556 x^6 - 68201 x^5 + 20860 x^4 - 4172 x^3 + 514 x^2 - 35 x + 1,\\
x ^{16} - 28 x ^{15} + 351 x ^{14} - 2604 x ^{13} + 12735 x ^{12} - 43295 x ^{11} + 105183 x ^{10} - 184976 x^9 + 236157 x^8 - 217812 x^7 + 143430 x^6 - 66161 x^5 + 20778 x^4 - 4257 x^3 + 532 x^2 - 36 x + 1,\\
x ^{16} - 29 x ^{15} + 376 x ^{14} - 2881 x ^{13} + 14533 x ^{12} - 50911 x ^{11} + 127385 x ^{10} - 230770 x^9 + 303849 x^8 - 289697 x^7 + 197940 x^6 - 95254 x^5 + 31441 x^4 - 6835 x^3 + 916 x^2 - 67 x + 2,\\
x ^{16} - 29 x ^{15} + 375 x ^{14} - 2857 x ^{13} + 14279 x ^{12} - 49348 x ^{11} + 121177 x ^{10} - 214063 x^9 + 272717 x^8 - 249314 x^7 + 161704 x^6 - 73113 x^5 + 22457 x^4 - 4506 x^3 + 554 x^2 - 37 x + 1,\\
x ^{16} - 29 x ^{15} + 375 x ^{14} - 2859 x ^{13} + 14320 x ^{12} - 49711 x ^{11} + 123001 x ^{10} - 219813 x^9 + 284591 x^8 - 265614 x^7 + 176517 x^6 - 81849 x^5 + 25677 x^4 - 5199 x^3 + 630 x^2 - 40 x + 1,\\
x ^{16} - 29 x ^{15} + 375 x ^{14} - 2857 x ^{13} + 14280 x ^{12} - 49366 x ^{11} + 121313 x ^{10} - 214625 x^9 + 274106 x^8 - 251427 x^7 + 163671 x^6 - 74197 x^5 + 22790 x^4 - 4557 x^3 + 557 x^2 - 37 x + 1,\\
x ^{16} - 28 x ^{15} + 350 x ^{14} - 2583 x ^{13} + 12543 x ^{12} - 42286 x ^{11} + 101792 x ^{10} - 177300 x^9 + 224146 x^8 - 204688 x^7 + 133439 x^6 - 60943 x^5 + 18971 x^4 - 3869 x^3 + 487 x^2 - 34 x + 1,\\
x ^{16} - 29 x ^{15} + 375 x ^{14} - 2858 x ^{13} + 14300 x ^{12} - 49539 x ^{11} + 122165 x ^{10} - 217274 x^9 + 279559 x^8 - 259009 x^7 + 170802 x^6 - 78665 x^5 + 24591 x^4 - 4995 x^3 + 614 x^2 - 40 x + 1,\\
x ^{17} - 31 x ^{16} + 433 x ^{15} - 3608 x ^{14} + 20017 x ^{13} - 78160 x ^{12} + 221435 x ^{11} - 462611 x ^{10} + 717469 x^9 - 825627 x^8 + 700247 x^7 - 432227 x^6 + 190457 x^5 - 58254 x^4 + 11862 x^3 - 1504 x^2 + 105 x - 3,\\
x ^{17} - 31 x ^{16} + 433 x ^{15} - 3607 x ^{14} + 19991 x ^{13} - 77864 x ^{12} + 219489 x ^{11} - 454402 x ^{10} + 694109 x^9 - 779715 x^8 + 637464 x^7 - 372774 x^6 + 152065 x^5 - 41797 x^4 + 7377 x^3 - 780 x^2 + 44 x - 1,\\
x ^{17} - 31 x ^{16} + 433 x ^{15} - 3608 x ^{14} + 20016 x ^{13} - 78141 x ^{12} + 221279 x ^{11} - 461883 x ^{10} + 715331 x^9 - 821490 x^8 + 694870 x^7 - 427515 x^6 + 187693 x^5 - 57190 x^4 + 11605 x^3 - 1469 x^2 + 103 x - 3,\\
x ^{17} - 30 x ^{16} + 404 x ^{15} - 3233 x ^{14} + 17157 x ^{13} - 63818 x ^{12} + 171512 x ^{11} - 338436 x ^{10} + 493528 x^9 - 531403 x^8 + 419445 x^7 - 239493 x^6 + 96981 x^5 - 27083 x^4 + 5010 x^3 - 577 x^2 + 37 x - 1,\\
x ^{18} - 31 x ^{17} + 434 x ^{16} - 3634 x ^{15} + 20320 x ^{14} - 80254 x ^{13} + 231007 x ^{12} - 493187 x ^{11} + 787795 x ^{10} - 943937 x^9 + 846712 x^8 - 565450 x^7 + 278719 x^6 - 100140 x^5 + 25720 x^4 - 4567 x^3 + 527 x^2 - 35 x + 1,\\
x ^{19} - 34 x ^{18} + 527 x ^{17} - 4935 x ^{16} + 31194 x ^{15} - 140864 x ^{14} + 469184 x ^{13} - 1173625 x ^{12} + 2224717 x ^{11} - 3203723 x ^{10} + 3495285 x^9 - 2868329 x^8 + 1750564 x^7 - 782324 x^6 + 250793 x^5 - 56098 x^4 + 8421 x^3 - 799 x^2 + 43 x - 1,\\
x ^{20} - 36 x ^{19} + 592 x ^{18} - 5901 x ^{17} + 39894 x ^{16} - 193914 x ^{15} + 700960 x ^{14} - 1922589 x ^{13} + 4046847 x ^{12} - 6570299 x ^{11} + 8228620 x ^{10} - 7915862 x^9 + 5800824 x^8 - 3198214 x^7 + 1304348 x^6 - 384755 x^5 + 79651 x^4 - 11093 x^3 + 974 x^2 - 48 x + 1,\\
x ^{21} - 38 x ^{20} + 662 x ^{19} - 7017 x ^{18} + 50652 x ^{17} - 264096 x ^{16} + 1029529 x ^{15} - 3064808 x ^{14} + 7056359 x ^{13} - 12651679 x ^{12} + 17706626 x ^{11} - 19317587 x ^{10} + 16350844 x^9 - 10650739 x^8 + 5276504 x^7 - 1955895 x^6 + 530450 x^5 - 102025 x^4 + 13307 x^3 - 1099 x^2 + 51 x - 1,\\
x ^{21} - 37 x ^{20} + 629 x ^{19} - 6520 x ^{18} + 46121 x ^{17} - 236127 x ^{16} + 905630 x ^{15} - 2657404 x ^{14} + 6041615 x ^{13} - 10714547 x ^{12} + 14856377 x ^{11} - 16082796 x ^{10} + 13529730 x^9 - 8775980 x^8 + 4340293 x^7 - 1611879 x^6 + 440266 x^5 - 85918 x^4 + 11484 x^3 - 984 x^2 + 48 x - 1,\\
x ^{21} - 37 x ^{20} + 629 x ^{19} - 6520 x ^{18} + 46121 x ^{17} - 236128 x ^{16} + 905654 x ^{15} - 2657659 x ^{14} + 6043199 x ^{13} - 10720947 x ^{12} + 14874086 x ^{11} - 16117217 x ^{10} + 13577179 x^9 - 8822270 x^8 + 4371857 x^7 - 1626575 x^6 + 444767 x^5 - 86772 x^4 + 11574 x^3 - 988 x^2 + 48 x - 1,\\
x ^{21} - 38 x ^{20} + 661 x ^{19} - 6986 x ^{18} + 50215 x ^{17} - 260384 x ^{16} + 1008299 x ^{15} - 2978253 x ^{14} + 6796491 x ^{13} - 12066219 x ^{12} + 16707028 x ^{11} - 18019925 x ^{10} + 15072989 x^9 - 9703686 x^8 + 4755571 x^7 - 1747760 x^6 + 471956 x^5 - 91010 x^4 + 12021 x^3 - 1018 x^2 + 49 x - 1,\\
x ^{21} - 37 x ^{20} + 629 x ^{19} - 6519 x ^{18} + 46093 x ^{17} - 235774 x ^{16} + 902977 x ^{15} - 2644137 x ^{14} + 5994942 x ^{13} - 10595414 x ^{12} + 14631953 x ^{11} - 15768232 x ^{10} + 13201104 x^9 - 8520973 x^8 + 4194537 x^7 - 1551385 x^6 + 422456 x^5 - 82341 x^4 + 11026 x^3 - 951 x^2 + 47 x - 1,\\
x ^{21} - 38 x ^{20} + 662 x ^{19} - 7016 x ^{18} + 50623 x ^{17} - 263716 x ^{16} + 1026549 x ^{15} - 3049190 x ^{14} + 6998496 x ^{13} - 12495319 x ^{12} + 17393108 x ^{11} - 18847465 x ^{10} + 15823311 x^9 - 10210292 x^8 + 5006216 x^7 - 1836383 x^6 + 493503 x^5 - 94394 x^4 + 12327 x^3 - 1030 x^2 + 49 x - 1,\\
x ^{21} - 37 x ^{20} + 628 x ^{19} - 6488 x ^{18} + 45656 x ^{17} - 232062 x ^{16} + 881750 x ^{15} - 2557644 x ^{14} + 5735638 x ^{13} - 10012934 x ^{12} + 13642521 x ^{11} - 14494118 x ^{10} + 11961174 x^9 - 7616692 x^8 + 3707186 x^7 - 1361262 x^6 + 370329 x^5 - 72717 x^4 + 9906 x^3 - 878 x^2 + 45 x - 1,\\
x ^{21} - 37 x ^{20} + 630 x ^{19} - 6550 x ^{18} + 46529 x ^{17} - 239462 x ^{16} + 923949 x ^{15} - 2729046 x ^{14} + 6247852 x ^{13} - 11160284 x ^{12} + 15587664 x ^{11} - 16997691 x ^{10} + 14401727 x^9 - 9405147 x^8 + 4679771 x^7 - 1746253 x^6 + 478206 x^5 - 93254 x^4 + 12396 x^3 - 1049 x^2 + 50 x - 1,\\
x ^{22} - 40 x ^{21} + 740 x ^{20} - 8407 x ^{19} + 65683 x ^{18} - 374549 x ^{17} + 1614800 x ^{16} - 5380628 x ^{15} + 14047901 x ^{14} - 28969248 x ^{13} + 47361479 x ^{12} - 61398902 x ^{11} + 62932078 x ^{10} - 50697588 x^9 + 31811123 x^8 - 15352578 x^7 + 5603004 x^6 - 1511320 x^5 + 291944 x^4 - 38597 x^3 + 3254 x^2 - 154 x + 3,\\
x ^{22} - 39 x ^{21} + 701 x ^{20} - 7710 x ^{19} + 58104 x ^{18} - 318403 x ^{17} + 1314115 x ^{16} - 4175106 x ^{15} + 10350675 x ^{14} - 20180685 x ^{13} + 31052363 x ^{12} - 37709766 x ^{11} + 36034271 x ^{10} - 26939661 x^9 + 15625755 x^8 - 6953429 x^7 + 2339557 x^6 - 583805 x^5 + 105193 x^4 - 13152 x^3 + 1069 x^2 - 50 x + 1,\\
x ^{22} - 40 x ^{21} + 737 x ^{20} - 8305 x ^{19} + 64097 x ^{18} - 359585 x ^{17} + 1518936 x ^{16} - 4938470 x ^{15} + 12528509 x ^{14} - 24997863 x ^{13} + 39368654 x ^{12} - 48935809 x ^{11} + 47854274 x ^{10} - 36582907 x^9 + 21654374 x^8 - 9794508 x^7 + 3325702 x^6 - 827755 x^5 + 146179 x^4 - 17489 x^3 + 1322 x^2 - 56 x + 1,\\
x ^{24} - 43 x ^{23} + 858 x ^{22} - 10556 x ^{21} + 89748 x ^{20} - 560151 x ^{19} + 2661377 x ^{18} - 9851449 x ^{17} + 28843603 x ^{16} - 67443670 x ^{15} + 126654006 x ^{14} - 191494024 x ^{13} + 233061937 x ^{12} - 227750986 x ^{11} + 177843081 x ^{10} - 110184505 x^9 + 53645540 x^8 - 20266836 x^7 + 5843592 x^6 - 1257832 x^5 + 196066 x^4 - 21177 x^3 + 1479 x^2 - 59 x + 1,\\
x ^{24} - 43 x ^{23} + 859 x ^{22} - 10592 x ^{21} + 90344 x ^{20} - 566173 x ^{19} + 2702939 x ^{18} - 10059252 x ^{17} + 29622628 x ^{16} - 69679771 x ^{15} + 131629715 x ^{14} - 200131135 x ^{13} + 244773712 x ^{12} - 240118139 x ^{11} + 187934748 x ^{10} - 116465181 x^9 + 56566618 x^8 - 21250343 x^7 + 6070888 x^6 - 1290222 x^5 + 198069 x^4 - 21067 x^3 + 1455 x^2 - 58 x + 1,\\
x ^{24} - 43 x ^{23} + 859 x ^{22} - 10591 x ^{21} + 90311 x ^{20} - 565675 x ^{19} + 2698377 x ^{18} - 10030848 x ^{17} + 29495071 x ^{16} - 69251786 x ^{15} + 130533011 x ^{14} - 197955518 x ^{13} + 241407113 x ^{12} - 236043600 x ^{11} + 184084671 x ^{10} - 113642568 x^9 + 54978625 x^8 - 20575825 x^7 + 5859381 x^6 - 1242735 x^5 + 190755 x^4 - 20343 x^3 + 1414 x^2 - 57 x + 1,\\
x ^{24} - 43 x ^{23} + 858 x ^{22} - 10555 x ^{21} + 89714 x ^{20} - 559620 x ^{19} + 2656321 x ^{18} - 9818597 x ^{17} + 28689109 x ^{16} - 66899428 x ^{15} + 125187554 x ^{14} - 188434347 x ^{13} + 228088005 x ^{12} - 221442351 x ^{11} + 171617735 x ^{10} - 105437276 x^9 + 50878424 x^8 - 19053046 x^7 + 5451455 x^6 - 1167270 x^5 + 181746 x^4 - 19728 x^3 + 1396 x^2 - 57 x + 1,$
\vfill \eject
\n $
x ^{24} - 43 x ^{23} + 859 x ^{22} - 10592 x ^{21} + 90343 x ^{20} - 566143 x ^{19} + 2702527 x ^{18} - 10055814 x ^{17} + 29603116 x ^{16} - 69599896 x ^{15} + 131385672 x ^{14} - 199563191 x ^{13} + 243755457 x ^{12} - 238705077 x ^{11} + 186418058 x ^{10} - 115213090 x^9 + 55779678 x^8 - 20879404 x^7 + 5942419 x^6 - 1258437 x^5 + 192673 x^4 - 20477 x^3 + 1418 x^2 - 57 x + 1,\\
x ^{24} - 43 x ^{23} + 858 x ^{22} - 10555 x ^{21} + 89714 x ^{20} - 559619 x ^{19} + 2656292 x ^{18} - 9818214 x ^{17} + 28686052 x ^{16} - 66882928 x ^{15} + 125123719 x ^{14} - 188251302 x ^{13} + 227691298 x ^{12} - 220785712 x ^{11} + 170784613 x ^{10} - 104628729 x^9 + 50282707 x^8 - 18724259 x^7 + 5318251 x^6 - 1128813 x^5 + 174174 x^4 - 18779 x^3 + 1329 x^2 - 55 x + 1,\\
x ^{24} - 43 x ^{23} + 858 x ^{22} - 10555 x ^{21} + 89714 x ^{20} - 559619 x ^{19} + 2656293 x ^{18} - 9818240 x ^{17} + 28686358 x ^{16} - 66885085 x ^{15} + 125133873 x ^{14} - 188285001 x ^{13} + 227772450 x ^{12} - 220929536 x ^{11} + 170972989 x ^{10} - 104810395 x^9 + 50410312 x^8 - 18788358 x^7 + 5340650 x^6 - 1134036 x^5 + 174934 x^4 - 18840 x^3 + 1331 x^2 - 55 x + 1,\\
x ^{24} - 43 x ^{23} + 859 x ^{22} - 10592 x ^{21} + 90344 x ^{20} - 566173 x ^{19} + 2702939 x ^{18} - 10059252 x ^{17} + 29622629 x ^{16} - 69679793 x ^{15} + 131629930 x ^{14} - 200132368 x ^{13} + 244778330 x ^{12} - 240130036 x ^{11} + 187956383 x ^{10} - 116493244 x^9 + 56592561 x^8 - 21267254 x^7 + 6078496 x^6 - 1292501 x^5 + 198497 x^4 - 21112 x^3 + 1457 x^2 - 58 x + 1,\\
x ^{24} - 43 x ^{23} + 858 x ^{22} - 10556 x ^{21} + 89748 x ^{20} - 560149 x ^{19} + 2661319 x ^{18} - 9850685 x ^{17} + 28837542 x ^{16} - 67411297 x ^{15} + 126530720 x ^{14} - 191148153 x ^{13} + 232333528 x ^{12} - 226587839 x ^{11} + 176429798 x ^{10} - 108880119 x^9 + 52737121 x^8 - 19795079 x^7 + 5664221 x^6 - 1209247 x^5 + 187074 x^4 - 20114 x^3 + 1408 x^2 - 57 x + 1,\\
x ^{26} - 47 x ^{25} + 1032 x ^{24} - 14073 x ^{23} + 133650 x ^{22} - 939584 x ^{21} + 5074160 x ^{20} - 21562027 x ^{19} + 73265558 x ^{18} - 201223170 x ^{17} + 449829027 x ^{16} - 821782544 x ^{15} + 1228894791 x ^{14} - 1503592842 x ^{13} + 1501724324 x ^{12} - 1219221943 x ^{11} + 799781540 x ^{10} - 420445551 x^9 + 175258436 x^8 - 57136292 x^7 + 14310194 x^6 - 2688900 x^5 + 366902 x^4 - 34680 x^3 + 2108 x^2 - 72 x + 1]$
\vskip1cm
\n
coef=[
0.53793383, 0.46790921, 0.06201384, 0.17528252, 0.00244359, 0.00594619, 0.06656685, 0.00415865, 0.00162666, 0.00026381, 0.02268826, 0.02087990, 0.00171650, 0.00670572,
0.00494538, 0.00435669, 0.00083999, 0.00026780, 0.00126922, 0.00059158, 0.00116447, 0.00026057, 0.00148560, 0.00074059, 0.00057768, 0.00004154, 0.00054812, 0.00222604,
0.00155412, 0.00207869, 0.00033049, 0.00022743, 0.00177254, 0.00018066, 0.00396194, 0.00036795, 0.00153950, 0.00006089, 0.00001082, 0.00024610, 0.00026106, 0.00013027,
0.00011303, 0.00118914, 0.00121424, 0.00072737, 0.00141169, 0.00133850, 0.00031424, 0.00009481, 0.00053529, 0.00022239, 0.00027012, 0.00265436, 0.00013343, 0.00008118,
0.00058771, 0.00110371, 0.00034140, 0.00057026, 0.00039390, 0.00034566, 0.00061844, 0.00167209, 0.00048967, 0.00053393, 0.00014761, 0.00026572, 0.00004545, 0.00051001,
0.00023718, 0.00053096, 0.00103947, 0.00040674, 0.00006127, 0.00017196, 0.00015398, 0.00015069, 0.00029189, 0.00059635, 0.00184900, 0.00027761, 0.00050559, 0.00005135,
0.00006235, 0.00015616, 0.00001694, 0.00041928, 0.00077785, 0.00096075, 0.00047160, 0.00097929, 0.00010843, 0.00006276, 0.00052579, 0.00007747, 0.00018519]

\m
\n
UMR CNRS 7502. IECL, Universit\'e de Lorraine, site de Metz, D\'epartement de Math\'ematiques, UFR MIM, Ile du Saulcy, CS 50128. 57045 METZ cedex 01. FRANCE\\
\n
E-mail address : valerie.flammang@univ-lorraine.fr

\end{document}